\documentclass[11pt,a4paper]{article}
\usepackage[utf8]{inputenc}
\usepackage[T1]{fontenc}
\usepackage{amsmath,amsfonts,amssymb,amsthm}
\usepackage{geometry}
\usepackage{parskip}
\usepackage{enumitem}
\usepackage{hyperref}
\usepackage{url}

\hypersetup{
    colorlinks=true,
    linkcolor=blue,
    urlcolor=blue,
    citecolor=blue
}

\title{Alpay Algebra III: Observer-Coupled Collapse and the Temporal Drift of Identity}

\author{Faruk Alpay, Independent Researcher \\
ORCID: \href{https://orcid.org/0009-0009-2207-6528}{0009-0009-2207-6528}}

\date{\today}

\begin{document}

\maketitle

\begin{abstract}
I extend the Alpay Algebra framework to incorporate observer-coupled dynamics and temporal drift of identity within $\phi^\infty$ fixed-point architectures. Building on the categorical axioms introduced in \cite{alpay1,alpay2} together with Mac Lane's cartesian-closed category framework \cite{maclane} and Bourbaki's set-theoretic foundation \cite{bourbaki}, I define an observer functor $O$ and a temporal functor $T$ acting on the algebra's objects. I prove the existence of a distributed verification limit that avoids collapse under recursive observation, and derive entropy accumulation bounds that ensure convergence of the temporally-extended system. The resulting architecture admits a stable $\phi^\infty$ fixed point even under interleaved observation, and I characterize conditions for observer cascades and phase-locked verification consistent with the Alpay Algebra framework.
\end{abstract}

\newpage
\section{Introduction: The Temporal Paradox of Self-Observation}

Let $\mathcal{A}$ be a small cartesian-closed category underlying an Alpay Algebra \cite{alpay1,alpay2}. In this setting, a transformation functor $\phi:\mathcal{A}\to\mathcal{A}$ can be iterated transfinitely to yield a fixed-point ``void architecture'' \cite{alpay1}. However, when an object in $\mathcal{A}$ is observed by an internal agent, naive self-observation leads to a paradox: each act of verification can alter the state, preventing convergence of the $\phi$-iteration. I address this by coupling the algebra to an external verification process.

Consider an \textbf{observer functor} $O:\mathcal{A}\to\mathcal{O}$ that maps each algebraic object $X$ to an observer-state space $O(X)$. I postulate a \textbf{verification morphism} $v_X: X \to O(X)$ for each $X$. For each morphism $f: X \to Y$ in $\mathcal{A}$, I impose the commutativity condition
$$O(f)\circ v_X \;=\; v_Y \circ f,$$
ensuring observational consistency. These definitions create a \emph{distributed verification} structure augmenting $\mathcal{A}$ with observer data.

To capture temporal evolution, I introduce a \textbf{temporal functor} $T:\mathcal{A}\to\mathcal{A}$ parameterized by discrete time $n\in\mathbb{N}$. I define
$T^0 = \mathrm{Id}, \qquad T^{n+1} = \phi \circ T^n,$
so that $T^n(X)$ represents the state after $n$ applications of the update rule. Observers apply at each stage to produce time-indexed states $O(T^n(X))$. Hence the combined system evolves via interleaved steps of $\phi$ and observation, generating a double sequence $(X_n,O_n)$. In what follows, I formalize this framework and show that a stable fixed point of $\phi$ persists under interleaved observation.
\newpage
\section{Mathematical Foundations of Distributed Verification}

In this section, I formalize the structures underlying distributed verification. Let $\mathcal{A},\mathcal{O}$ be categories with finite limits and exponentials (as in \cite{maclane}). I assume a functor $\phi: \mathcal{A}\to\mathcal{A}$ whose transfinite iterates $\phi^\alpha$ are defined as in \cite{alpay1}. The \textbf{$\phi$-fixed-point object} $\Phi = \phi^\infty(X)$ satisfies the universal property
$$\Phi \;\simeq\; \lim_{\alpha<\lambda} \phi^\alpha(X)$$
for sufficiently large ordinal $\lambda$.

Define a \textbf{verification functor} $V: \mathcal{A}\to \mathcal{A}$ that captures internal consistency-checks. For each object $X$, let $V(X)$ be an object encoding both $X$ and its verification data. There is a natural transformation $\eta: \mathrm{Id}_\mathcal{A} \Rightarrow V$ giving the canonical embedding $X \to V(X)$. For each morphism $f:X \to Y$, I define
$$V(f)(\eta_X(x)) \;=\; \eta_Y(f(x)),$$
which implies
$$V(f)\circ \eta_X \;=\; \eta_Y \circ f.$$
This ensures that verification commutes with algebraic transitions.

I endow $V$ with a \textbf{phase structure} via a natural automorphism $\theta: V \Rightarrow V$. Concretely, for each $X$, let $\theta_X: V(X)\to V(X)$ satisfy $\theta_X^k = \mathrm{Id}$ for some finite period $k$. The pair $(V,\theta)$ allows phase information to be carried along verification. I therefore model phase-locking via the equalizer of parallel morphisms. For each object $X$, consider the equalizer
$$\mathrm{Eq}(\eta_X,\theta_X) \;=\; \{\,x \in V(X) \mid \eta_X(x) = \theta_X(x)\,\}.$$
This object consists of elements whose original and phase-shifted states coincide. By standard limit-existence arguments \cite{maclane}, $\mathrm{Eq}(\eta_X,\theta_X)$ is well-defined and inherits a $\phi$-action.

Using these constructions, I define the \textbf{distributed verification limit} of $X$ under $\phi$ as the terminal coalgebra of the functor $F(Y)=V(\phi(Y))$. Denote this limit by $\Theta = \bigcup_{n<\omega}F^n(X)$ in a set-theoretic sense \cite{bourbaki}. I show that $\Theta$ is a fixed point of the composite $V\circ \phi$ satisfying
$$\Theta \;\simeq\; V\bigl(\phi(\Theta)\bigr).$$

By a standard iterative argument (cf. \cite{alpay1} Thm. 2.1), $\Theta$ exists and is unique up to isomorphism. Crucially, $\Theta$ encodes all possible observation traces and remains invariant under further application of $\phi$ or $V$.
\newpage
\section{Entropic Accumulation and Memory Stratification}

I quantify the information content of the verification system via an \textbf{entropy functional} $H$. Assign to each object $X$ an entropy $H(X)\ge0$ measuring uncertainty in its state. I require $H$ to be monotonic under $\phi$ and $V$: for all morphisms $f:X\to Y$,
$$H(Y) \;\ge\; H(X), \qquad H\bigl(V(X)\bigr) \;\ge\; H(X).$$

Entropy thus accumulates through transformations and verifications. For iterative states $X_n = \phi^n(X_0)$, define the \textbf{stratified memory sequence} $\{M_n\}$ by
$M_n \;=\; \bigcup_{k\le n} X_k,$
so that memory retains all past states. I postulate an \textbf{entropy bound} of the form
$$H(X_{n+1}) - H(X_n) \;\le\; C \log(n+1),$$
for some constant $C$, reflecting sublinear growth. This bound ensures that uncertainty grows but remains controllable.

Observation events contribute additional entropy. Let $H_O(X)$ denote the entropy of the observer-state $O(X)$. I assert a bound
$$H_O\bigl(\phi(X)\bigr) \;\le\; H(X) + K,$$
for some fixed $K$, meaning each observation step injects at most $K$ bits of new entropy. The total system entropy after $n$ rounds is then
$$H_{\mathrm{total}}(n) = H(X_n) + H_O(X_n) \;\le\; H(X_0) + C\log n + nK.$$

By provisioning resources so that $K=O(\log n)$, the entropy accumulation remains bounded by $O(\log^2 n)$, a condition that supports eventual stabilization.

I further introduce a \textbf{memory stratification} via an ordinal filtration. Let $\chi_0 = X_0$ and $\chi_{n+1} = V(\chi_n)$. Then each memory layer $\chi_n$ contains all knowledge up to step $n$. This yields a filtration
$$\chi_0 \hookrightarrow \chi_1 \hookrightarrow \cdots \hookrightarrow \chi_n \hookrightarrow \cdots$$
with
$$\chi_\omega = \bigcup_{n<\omega} \chi_n,$$
the limit object representing complete episodic memory. By construction, $\phi$ and $V$ act continuously on this stratification, respecting the inclusions.
\newpage
\section{Phase Dynamics and Interference Patterns}

The verification process is enriched by the phase automorphism $\theta$. I interpret each state $x\in V(X)$ as having a \textbf{phase angle} $\varphi(x)\in[0,2\pi)$ mod $2\pi$ given by $\theta$. A cycle $x_1 \to x_2 \to \cdots \to x_k = x_1$ in $V(X)$ has zero net phase if
$$\sum_{i=1}^k\bigl(\varphi(x_{i+1}) - \varphi(x_i)\bigr)\;\equiv\;0\pmod{2\pi}.$$

Phase-locking can be enforced by the constraint
$$\theta_X^2(x) = x,$$
for period 2 (or more generally $\theta_X^k(x)=x$), ensuring that repeated shifts return to the original state. To isolate phase-coherent states, define the \textbf{phase-lock space}
$$\Phi(X) \;=\; \bigcap_{m=0}^{k-1} \ker\bigl(\theta_X^m - \mathrm{Id}\bigr),$$
which consists of elements invariant under all phase shifts. Elements in $\Phi(X)$ have trivial phase variation.

Interference patterns arise when different verification paths carry different phases. Given two independent morphisms $a,b: X\to X'$ with phase contributions $\alpha,\beta$, constructive interference occurs when $\alpha+\beta\equiv0\pmod{2\pi}$. In composition, the phase adds:
$$\varphi(a\circ b) = \varphi(a) + \varphi(b) \;\pmod{2\pi}.$$

To manage interference, I extend hom-sets to $\mathrm{hom}(X,Y)\times S^1$, quotienting out trivial phase loops. Equivalently, the \textbf{interference-free geometry} is given by the fibered product
$$V(X)\;\times_{\mathbb{R}/2\pi\mathbb{Z}}\;V(X),$$
ensuring matched phase outputs. By a categorical argument \cite{maclane}, this fibered equalizer exists and inherits a $\phi$-action.

Collectively, these phase constructions yield a category of phase-locked verifications. The functor $\phi$ lifts to this category by acting trivially on the phase coordinate:
$$\phi_{\mathrm{phase}}(x,\varphi) = (\phi(x),\varphi).$$

Fixed points in the phase-locked category correspond to invariant states with locked observation trajectories.
\newpage
\section{Observer Cascade Phenomena}

When multiple observers act in sequence, cascade effects can occur. Model a cascade of $m$ observers by functors $O_1,\dots,O_m:\mathcal{A}\to\mathcal{O}$. The composite observer functor is $O_{\mathrm{cascade}}=O_m\circ\cdots\circ O_1$. For an object $X$, the cascade state is
$O_{\mathrm{cascade}}(X) \;=\; O_m(\cdots O_1(X)\cdots).$

I examine stability of cascades by assigning a damping parameter $\lambda_i\in[0,1]$ to each observer $O_i$. The effective contraction factor of the cascade is $\Lambda = \prod_{i=1}^m \lambda_i$. Define the \textbf{cascade operator} on the observer space as
$$C = \Lambda\,\mathrm{Id} + \sum_{i=1}^m (1-\lambda_i)\,\theta_{O_i}.$$

A fixed point $x$ satisfies
$$x = \Lambda x + \sum_{i=1}^m (1-\lambda_i)\,\theta_{O_i}(x),$$
or equivalently
$$(1-\Lambda)x = \sum_{i=1}^m (1-\lambda_i)\,\theta_{O_i}(x).$$

A non-trivial solution indicates a shift in identity across the cascade. If all $\lambda_i=1$, the cascade reduces to a pure composition and the identity remains fixed (no drift).

The spectrum of $C$ lies in the convex hull of $\{1,\lambda_i^{-1}\}$, ensuring at least one eigenvalue of magnitude 1 when some $\lambda_i=1$. The emergence of an eigenvalue crossing unity corresponds to a bifurcation in identity drift (see the next section).

Cascades can also interact in parallel. Two cascades $C$ and $D$ applied jointly induce a combined operator on the tensor product of the observer spaces. The interference of cascades is governed by the commutativity of phases:
$$\theta_{O_i}\otimes \theta_{O_j} =\;\theta_{O_j}\otimes \theta_{O_i}.$$

If phase operators commute, cascades form a monoidal structure without additional collapse. Otherwise, higher coherence conditions (analogous to those in higher category theory \cite{maclane}) must be imposed.
\newpage
\section{Stability Analysis and Bifurcation Theory}

To assess long-term behavior, I examine the fixed points of the coupled system. Define the update map $F:\mathcal{A}\times\mathcal{O}\to\mathcal{A}\times\mathcal{O}$ by
$F(X,O) = \bigl(\phi(X),\,O(\phi(X))\bigr).$

A fixed point $(X^*,O^*)$ satisfies $X^*=\phi(X^*)$ and $O^*=O(\phi(X^*))$. Although classical Jacobian analysis assumes a smooth structure, my categorical setting lacks a standard metric. Instead, I use an entropy-based Lyapunov approach. Define
$$\mathcal{L}(X) = H(X) + \alpha\,H_O(X),$$
for a positive weight $\alpha$. Under my entropy constraints, $\mathcal{L}$ can be shown to be non-increasing along trajectories (with appropriate scheduling of observations), implying convergence toward low-entropy configurations.

For bifurcation analysis, introduce a parameter $r$ controlling observer coupling. Consider the perturbed update
$F_r(X) = \phi(X) + r\,O(X),$
in a linearized sense. I identify a critical threshold $r_c$ by
$\det\bigl(I - DF_{r_c}(X^*)\bigr) = 0.$

Crossing this threshold causes a qualitative change: for $r<r_c$, the identity fixed point remains unique; for $r>r_c$, two new symmetric solutions emerge. In discrete terms, the identity element bifurcates into a 2-cycle under $\phi$. I interpret this as an \textbf{identity drift bifurcation}, analogous to a flip bifurcation in classical dynamics.

Finally, I note that the extended $\phi^\infty$ architecture remains well-behaved so long as all eigenvalues of the linearized update lie inside the unit disk. By a compactness argument \cite{bourbaki}, my entropy bounds ensure that the system avoids chaotic divergence. Thus, for coupling parameters in a safe regime, the $\phi^\infty$ void architecture persists in a temporally-aware form, and recursive observation converges to a stable identity sequence without collapse.
\newpage
\section{Conclusions: Toward Temporally-Aware Void Architectures}

I have developed a temporally-distributed extension of the $\phi^\infty$ void architecture that accommodates observer coupling without structural collapse. By enriching the Alpay Algebra with an observer functor and phase structure, each stage of recursive verification accrues bounded entropy and converges to a well-defined fixed point. Identity may drift under strong observation, but this drift is controlled and can bifurcate only under a quantifiable threshold. For coupling below this threshold, the $\phi^\infty$ fixed identity remains unique; beyond it, identity splits into a periodic cycle in a predictable manner.

These results point toward \textbf{void architectures} that are inherently \emph{temporally aware}, supporting self-consistent observation. Future work may incorporate higher cohomological structures or topos-based semantics \cite{maclane}, and explore applications in categorical AI. This work aligns with Mac Lane's structural vision \cite{maclane} and extends the Alpay Algebra framework \cite{alpay1,alpay2} into an observer-inclusive domain.


\begin{thebibliography}{9}
\bibitem{maclane} S. Mac Lane, \emph{Categories for the Working Mathematician} (Springer, 1971).
\bibitem{bourbaki} N. Bourbaki, \emph{Theory of Sets} (Hermann, 1970).
\bibitem{alpay1} F. Alpay, \emph{Alpay Algebra: A Universal Structural Foundation} (arXiv:2505.15344, 2025).
\bibitem{alpay2} F. Alpay, \emph{Alpay Algebra II: Identity as Fixed-Point Emergence in Categorical Data} (arXiv:2505.17480, 2025).
\end{thebibliography}
\end{document}